\documentclass[a4paper,12pt]{article} 
\usepackage{amsmath, amsthm, amssymb}
\usepackage{url}

\usepackage{physics}
\usepackage{graphicx,lipsum}
\graphicspath{ {./downloads/} }

\usepackage[colorlinks,citecolor=red,urlcolor=blue,bookmarks=false,hypertexnames=true]{hyperref}
\usepackage{tikz}
\usetikzlibrary{calc}
\usetikzlibrary{shapes}
\usepackage[autostyle]{csquotes}
\makeatletter

\usepackage[colorlinks]{hyperref}
\usepackage[nameinlink,capitalize]{cleveref}
\newtheorem{theorem}{Theorem}[section]
\newtheorem{corollary}[theorem]{Corollary}
\newtheorem{lemma}[theorem]{Lemma}
 \newtheorem{proposition}[theorem]{Proposition}
 \newtheorem{remark}[theorem]{Remark}
\newtheoremstyle{named}{}{}{\itshape}{}{\bfseries}{.}{.5em}{\thmnote{#3's }#1}
\theoremstyle{named}

\theoremstyle{definition}
\newtheorem{definition}[theorem]{Definition}
\newtheorem{example}[theorem]{Example}

\usepackage{xspace}
\usepackage[margin=1.0in]{geometry}

\usepackage{titlesec}

\usepackage{mathtools}

\DeclarePairedDelimiterX{\inp}[2]{\langle}{\rangle}{#1, #2}
\titleformat{\chapter}
  {\Large\bfseries} 
  {}                
  {0pt}            
  {\huge}

\begin{document}

\begin{center}
\fontsize{13pt}{10pt}\selectfont
    \textsc{\textbf{GENERALIZATIONS OF GRADED PRIME IDEALS OVER GRADED NEAR RINGS}}
    \end{center}
\vspace{0.1cm}
\begin{center}
   \fontsize{12pt}{10pt}\selectfont
    \textsc{{\footnotesize  Malik Bataineh*, Tamem Al-shorman and Eman Al-Kilany }}
\end{center}
\vspace{0.2cm}

\begin{abstract}
  This paper considers graded near-rings over a monoid G as a generalizations of the graded rings over groups, introduce certain innovative graded weakly prime ideals and graded almost prime ideals as a generalizations of graded prime ideals over graded near-rings, and explore their various properties and their generalizations in graded near-rings.    
\end{abstract}

\section{INTRODUCTION}
Throughout this article, G will be an abelian group  with identity e and R be a commutative ring with nonzero unity 1 element. R is called a G-graded ring if $ R= \bigoplus\limits_{g \in G} R_g$   with the property $R_gR_h\subseteq R_{gh}$ for all $g,h \in G$, where $R_g$ is an additive subgroup of R for all $g\in G$. The elements of $R_g$ are called homogeneous of degree g. If $x\in R$, then $x$ can be written uniquely as $\sum\limits_{g\in G} x_g$, where $x_g$ is the component of $x$ in $R_g$. The set of all homogeneous elements of R is $h(R)= \bigcup\limits_{g\in G} R_g$. Let P be an ideal of a G-graded ring R. Then P is called a graded ideal if $P=\bigoplus\limits_{g\in G}P_g$, i.e, for $x\in P$ and  $x=\sum\limits_{g\in G} x_g$ where $x_g \in P_g$ for all $g\in G$. An ideal of a G-graded ring is not necessary graded ideal (see \cite{abu2019graded}). The concept of graded prime ideals and its generalizations have an indispensable role in commutative G-graded rings. 

Near-rings are generalizations of rings in which addition is not necessarily abelian and only one distributive law holds. They arise in a natural way in the study of mappings on groups: the set M(G) of all maps of a group (G; +) into itself endowed with point-wise addition and composition of functions is a near-ring. For general background on the theory of near-rings, the monographs written by Pilz \cite{pilz2011near} and Meldrum \cite{meldrum1985near} should be referred. The definition of a near-ring ($N$, +, $\times$) is a set $N$ with two binary + and $\times$ that satisfy the following axioms:
\\
(1) ($N$, +) is a group.
\\
(2) ($N$, $\times$) is semi group. (semi group: a set together with an associative binary operation).
\\(3) $\times$ is right distributive over + (i.e. $(a+b)\times y = ay + by$).

The graded rings were introduced by Yoshida in \cite{yoshida1955homogeneous}. Also, graded near-rings were introduced and studied by Dumitru, Nastasescu and Toader in \cite{dumitru2016graded}. Let G be a multiplicatively monoid (an algebraic structure with a single associative binary operation) with identity. A near-ring $N$ is called a G-graded near-ring if there exists a family of additive normal subgroups $\{N_\sigma\}$ of $N$ satisfying that:
\\
(1) $N = \bigoplus\limits_{\sigma \in G} N_\sigma$.
\\
(2) $N_\sigma N_\tau \subseteq N_{\sigma \tau}$ for all $\sigma, \tau \in G$.
\\
A graded ideal P of a G-graded ring R is said to be graded prime ideal of R if $ab \in P$, where $a, b \in h(R)$, then $a \in P$ or $b \in P$. Graded prime ideals have been generalized to graded weakly prime ideals and graded almost prime ideals. In \cite{atani2006graded}, a graded ideal P of R is said to be graded weakly prime ideal of R if $0\neq ab \in P$, where $a, b \in h(R)$, then $a\in P$ or $b \in P$. We say that a graded ideal P of R is a graded almost prime ideal of R if $ab \in P-[P^2 \cap R]$, where $a, b \in h(R)$, then $a\in P$ or $b \in P$ (see \cite{jaber2008almost}).

Bataineh, Al-Shorman and  Al-Kilany in \cite{bataineh2022graded},  defined the concept of graded prime ideals over graded near-rings. A graded ideal P of a graded near-ring $N$ is said to be a graded prime ideal of $N$ if whenever $IJ \subseteq P$, then either $I \subseteq P$ or $J \subseteq P$, for any graded ideals $I$ and $J$ in $N$. In Section Tow, we introduced the concept of graded weakly prime ideals in graded near-rings. We say that P is a graded weakly prime ideal of $N$ if whenever $\{0\} \neq IJ \subseteq P$, then either $I \subseteq P$ or $J \subseteq P$, for any graded ideals $I$ and $J$ in $N$. In Section Three, we introduce the concept of graded almost prime ideals in graded near-ring. We say that P is a graded almost prime ideal of $N$ if whenever $IJ\subseteq P$ and $IJ \not\subseteq (P^2 \cap N)$, then either $I \subseteq P$ or $J \subseteq P$, for any graded ideals $I$ and $J$ in $N$.

\section{ GRADED WEAKLY PRIME IDEALS OVER GRADED NEAR RINGS}

In this section, we introduce graded weakly prime ideals graded over near-rings concept and study their basic properties.

\begin{definition}\label{def1}
Let G be a multiplucative monoid group with identity element and N be a G-graded near-ring. A graded ideal P of N is called graded weakly prime ideal of N if whenever $\{0\} \neq IJ \subseteq P$, then either $I \subseteq P$ or $J \subseteq P$, for any graded ideals I and J in N.
\end{definition}

\begin{example}\label{ex1}
Consider the ring ($\mathbf{Z}_{12}$, +, $\times$) is a near-ring with $G=\{0,1\}$ is a group under (+), where (+) defined as 0+0=0, 0+1=1, 1+0=1 and 1+1=1.  
\\
Let N be a G-graded near-ring defined by $N_0= \mathbf{Z}_{12}$ and $N_1 =\{0\}$. Note that the graded ideals $P_1= \{0\}$, $P_2=\{0,2,4,6,8,10\}$ and $P_3=\{0,3,6,9\}$ are graded weakly prime ideals of N.
\end{example}

\begin{remark}\label{rem1}
Every graded prime ideals over graded near-rings is a graded weakly prime ideals over graded near-rings. However, the converse is not true. For Example \ref{ex1}, $P_1$ is a graded weakly prime ideal of N but not graded prime ideal of N.
\end{remark}

The following theorem and corollary state that graded weakly prime ideals of N are graded prime ideals of N when certain conditions are met.

\begin{theorem}\label{thm1}
Let N be a G-graded near-ring and P be a graded weakly prime ideal of N. If P is not graded prime ideal of N, then $P^2 \cap N = \{0\}$.
\\
\\
\textbf{Proof.} Suppose that $P^2 \cap N \neq \{0\}$. It is observed that P is graded prime ideal of N. Let I and J be a graded ideals of N such that $IJ \subseteq P$. If $IJ \neq \{0\}$, then $I \subseteq P$ or $J\subseteq P$ since P is a graded weakly prime ideal of N. So, it could be assumed that $IJ=\{0\}$. Since  $P^2 \cap N \neq \{0\}$, so there exists $p,q \in P$ such that $<p><q> \neq 0$ and so $(I+<p>)(J+<q>) \neq \{0\} $. Suppose that $(I+<p>)(J+<q>) \not\subseteq P$, then there exists $i \in I$, $j\in J$, $p_0\in <p>$ and $q_0 \in <q>$ such that $(i+p_0)(j+q_0) \not\in P$ which implies that $i(j+q_0) \not\in P$, but $i(j+q_0)=i(j+q_0) - ij \in P$ since $IJ = \{0\}$. This is a contradiction. Thus, $\{0\} \neq (I+<p>)(J+<q>) \subseteq P $ which implies that $I \subseteq P$ or $J \subseteq P$. 
\end{theorem}

\begin{corollary}\label{coro1}
Let N be a G-graded near-ring and let P be a graded ideal of N such that $P^2 \cap N \neq \{0\}$. Then P is graded prime ideal of N if and only if P is graded weakly prime ideal of N.
\\
\\
\textbf{Proof.} Let P be a graded ideal of N such that $P^2 \cap N \neq  \{0\}$. By Theorem \ref{thm1}, if P is a graded weakly prime ideal of N, than P is a graded prime ideal of N. Also, by Remark \ref{rem1}, if P is a graded prime ideal of N, then P is a graded weakly prime ideal of N.
\end{corollary}

\begin{remark}\label{rem2}
It is not necessary that P is graded weakly prime ideal of N such that $P^2 \cap N = \{0\}$. Let N be a graded near-ring which is defined in Example \ref{ex1} and let $P=\{0,6\}$, Note that $P^2 \cap N =\{0\}$, but P is not graded weakly prime ideal of N.
\end{remark}

The next proposition gives an interesting case where graded weakly prime ideals lead to graded prime ideals in a graded near-ring.

\begin{proposition}\label{prop1}
Let N be a G-graded near-ring and P be a graded ideal of N. If P is a graded weakly prime ideal of N and $(\{0\} : P ) \subseteq P$, then P is a grade prime ideal of N.
\\
\\
\textbf{Proof.} Suppose that P is not a graded prime ideal of N, then there exists $I \not\subseteq P$ and $J \not\subseteq P$ satisfying that $IJ \subseteq P$, where I and J are two graded ideals of N. If $IJ \neq \{0\}$, then it is completed. So, it is assume that $IJ = \{0\}$. Note that $IJ \subseteq P$ since if an element belongs to $IP$, then it belongs to both N and P. Consider $I(J+P) \subseteq P$ if $I(J+P) \neq \{0\}$, then either $I \subseteq P$ or $J\subseteq P$, this is a contradiction. Otherwise, $I(J+P) = \{0\}$, then $IP = \{0\}$ implies $I \subseteq (\{0\} :P) \subseteq P$. 
\end{proposition}

\begin{theorem}\label{thm2}
Let N be a G-graded near-ring and P be a graded weakly prime ideal of N. If $IJ = \{0\}$ with $I \not\subseteq P$ and $J \not\subseteq P$ where I and J are two graded ideals of N, then $IP = PJ$.
\\
\\
\textbf{Proof.} Suppose that there exists $p \in P$ and $i \in I$ such that $ip = 0$. Then $\{0\} \neq I(J+<p>) \subseteq P$. But $I \not\subseteq P$ and $J+<p> \not\subseteq P$, which contradicts that P being graded weakly prime ideal of N. 
\end{theorem}

\begin{lemma}\label{lem1}
Let N be a G-graded near-ring. If P, I and J are graded ideals of N such that $P = I \cup J$, then P equals I or J.
\\
\\
\textbf{Proof.} Suppose P does not equal I nor J. Let $x \in P$ such that $x\in I$ but $x\not\in J$ and $y \in P$ such that $y\in J$ but 
$y\not\in I$. Since P is a graded ideal of N, $x-y\in P$. This implies $x-y \in I$ or $x-y \in J$. If $x-y\in I$, then $y \in I$ since I is a graded ideal of N, this is a contradiction. If $x-y\in J$, then $x\in J$ since J is a graded ideal of N, this is a contradiction. Therefore, P equals either I or J.
\end{lemma}

\begin{proposition}\label{prop2}
Let N be a G-graded near-ring and P be a graded ideal of N. Then the following are equivalent:
\\
(1) For $x,y \ and\ z \in N$ with $0\neq x(<y>+<z>) \subseteq P$, $x\in P$ or $y \ and\ z \in P$.
\\
(2) For $x\in N$ but $x \not \in P$ we have $(P:<x>+<y>) = P \cup (0:<x>+<y>)$ for any $y \in N$.
\\
(3) For $x\in N$ but $x \not\in P$ we have $(P:<x>+<y>) = P$ or $(P:<x>+<y>) = (0:<x>+<y>)$ for any $y\in N$.
\\
(4) P is a graded weakly prime ideal of N.
\\
\\
\textbf{Proof.} $(1) \Rightarrow (2)$: Let $t \in N$ and $t \in (P :< x > + < y >)$ for any y and x belongs to N but $x \not\in P$. Then $t(< x > + < y >) \subseteq P$. If $t(< x > + < y >) = 0$. Then $t \in (0 :< x > + < y >)$. Otherwise $0 \neq t(< x > + < y >) \subseteq P$. Thus, $t \in P$ by hypothesis. 
\\
$(2) \Rightarrow (3)$: It is following directly from Lemma \ref{lem1}.
\\
$(3) \Rightarrow (4)$: Let I and J be a graded ideal of N such that $IJ \subseteq P$. Suppose that $I \not\subseteq P$ and $J \not\subseteq P$. Then there exist $j \in J$ with $j \not\in P$. Now, it is claimed that $IJ = \{0\}$. Let $j_1 \in J$, then $I(<j>+<j_1> \subseteq P$, which implies $I \subseteq (P:<j>+<j_1>)$. Then by assumption, $I(<j>+<j_1>) =0 $ which gives $Ij_1 = \{0\}$. Thus $IJ = \{0\}$ and hence P is a graded weakly prime ideal of N.
\\
$(4) \Rightarrow (1)$: If $0 \neq x(<y>+<z>) \subseteq P$, then $\{0\} \neq <x>(<y>+<z>) \subseteq P$. Since P is a graded weakly prime ideal of N, there is $<x> \subseteq P$ or $<y>+<z> \subseteq P$. By assumption, $x, y \ and\ z \in N$. Hence $x\in P$ or $y \ and\ z \in P$. 
\end{proposition}

\begin{theorem}\label{thm3}
Let N be a G-graded near-ring and P be a graded ideal of N. Then the following are equivalent:
\\
(1) P is a graded weakly prime ideal of N.
\\
(2) For any ideals I and J in N with $P \subset I$ and $P \subset J$, then there is either $IJ=\{0\}$ or $IJ \not\subseteq P$.
\\
(3) For any ideals I and J in N with $I \not\subseteq P$ and $J \not\subseteq P$, then there is either $IJ=\{0\}$ or $IJ \not\subseteq P$.
\\
\\
\textbf{Proof.} $(1) \Rightarrow (2)$: Let I an J be two graded ideals of N with $P \subset I$, $P \subset J$ and $IJ\neq \{0\}$. Take $i \in I$ and $j \in J$ with $i \not \in P$ and $j \not \in P$, which implies that $\{0\} \neq <i><j> \not\subseteq P$ and hence $\{0\} \neq IJ \not\subseteq P$.
\\
$(2) \Rightarrow (3)$: Let I and J be a graded ideals of N with $I \not\subseteq P$ and $J \not\subseteq P$. Then there exists $i_1 \in I$ and $j_1 \in J$ such that $i_1 \not \in P$ and $j_1 \not\in P$. Suppose that $<i><j> \neq \{0\}$ for some $i \in I$ and $j \in J$. Then $(P+<i>+<i_1>)(P+<j>+<j_1>) \neq \{0\}$ and $P \subset (P+<i>+<i_1>) $ and $ P \subset (P+<j>+<j_1>)$. By hypothesis, $(P+<i>+<i_1>)(P+<j>+<j_1>) \not\subseteq P$. So, $<i>(P+<j>+<j_1>) +<i_1>(P+<j>+<j_1>) \not\subseteq P$. Hence there exists $i' \in <i>$, $i_1'\in <i_1>$, $j', j'' \in <j>$, $j_1', j_1'' \in <j_1>$ and $p_1, p_2 \in P$ such that $i'(p_1 + j' + j_1') + i_1'(p_2 + j'' + j_1'') \not\in P$. Thus  $i'(p_1 + j' + j_1') - i'(j' + j_1') + i'(j' + j_1') + i_1'(p_2 + j'' + j_1'') - i_1'(j'' + j_1'') + i_1'(j'' + j_1'') \not\in P$. But $i'(p_1 + j' + j_1') - i'(j' + j_1') \in P$ and $i_1'(p_2 + j'' + j_1'') - i_1'(j'' + j_1'') \in P$. Which implies neither $i'(j' + j_1')$ nor $i_1'(j'' + j_1'')$ belongs to P. Therefore, $IJ \not\subseteq P$.
\\
$(3) \Rightarrow (1) $: Follows directly from the definition of graded weakly prime ideals of N.
\end{theorem}

\begin{proposition}\label{prop3}
Let N be a G-graded near-ring, A be a totally ordered set and $(P_a)_{a\in A}$ be a family of graded weakly prime ideals of N with $P_a \subseteq P_b$ for any $a,b \in A$ with $a\leq b$. Then $P = \bigcap\limits_{a\in A} P_a$ is a graded weakly prime ideal of N.
\\
\\
\textbf{Proof.} Let I and J be two graded ideals of N with $\{0\} \neq IJ \subseteq P$, which implies for all $a \in A$ there is $IJ \subseteq P_a$. If there exists $a \in A$ such that $I \not\subseteq P_a$, then $J \subseteq P_a$. Hence for all $a \leq b$ there is $J \subseteq P_b$. If there exists $c < a$ such that $j \not \subseteq P_c$, then $I \subseteq P_c $ and then $I \subseteq P_a$, this is a contradiction. Hence for any $a \in A$, there is $J \subseteq P_a$. Therefore, $J \subseteq P$. 
\end{proposition}

\begin{proposition}\label{prop4}
Let N be a G-graded near-ring and P be an intersection of some graded weakly prime ideals of N. Then for any graded ideal I of N satisfying that $\{0\} \neq I^2 \subseteq P$ there is $I \subseteq P$.
\\
\\
\textbf{Proof.} Let $P_a$ be a set of graded weakly prime ideals of N, P be the intersection of $P_a$ and I be a graded ideal of N such that $\{0\} \neq I^2 \subseteq P$. Then $I^2$ is subset of each $P_a$ since $P_a$ is graded weakly prime ideal of N there is $I \subseteq P_a$. Therefore, $I \subseteq P$.
\end{proposition}

Next Example and Theorem \ref{thm4}, shows that the pre-image of a surjective homomorphism map of graded weakly prime ideal of N is not necessary to be graded weakly prime ideal of N, while the image of a surjective homomorphism map of graded weakly prime ideal of N which contains the kernal is graded weakly prime ideal of N.

\begin{example}\label{ex2}
Let G be the multiplicatively monoid which defined in Example \ref{ex1} and $N = \mathbf{Z}_8$ and $M=\mathbf{Z}_4$ be two G-graded near-rings where $N_0 = \mathbf{Z}_8$, $N_1 = \{0\}$, $M_0 = \mathbf{Z}_4$ and $M_1 = \{0\}$. Consider $\phi : N \rightarrow M$ where $\phi(x) =X$ is surjective homomorphism map. However, $\{0\}$ is a graded weakly prime ideal in M although $\phi^{-1}(\{0\}) = \{0, 4\}$ is not graded weakly prime ideal of N. 
\end{example}

\begin{lemma}\label{lem2}
Let N and M be two G-graded near-rings and $\phi$ be a surjective homomorphism from N into M. For any two graded ideals I and J of N if $IJ \neq \{0\}$, then $\phi^{-1}(I) \phi^{-1}(J) \neq \{0\}$.
\\
\\
\textbf{Proof.} Let I and J be two graded ideals of N such that $IJ \neq \{0\}$. Suppose that $\phi^{-1}(I) \phi^{-1}(J) = \{0\}$, then $\phi^{-1}(I) \phi^{-1}(J) = \phi^{-1}(IJ) = \{0\}$. Therefore, $\phi(\{0\}) = IJ$ which contradict with the fact that the image of zero is zero for any homomorphism map. Hence $\phi^{-1}(I) \phi^{-1}(J) \neq \{0\}$.
\end{lemma}

\begin{theorem}\label{thm4}
Let N and M be two G-graded near-rings and $\phi$ be a surjective homomorphism from N into M. Then the image of graded weakly prime ideal of N which contains the kernal of $\phi$ is a graded weakly prime ideal of M.
\\
\\
\textbf{Proof.} Suppose that $\{0\} \neq IJ \subseteq \phi (P)$ where I and J are graded ideals of M and P is a graded weakly prime ideals of N. By Lemme \ref{lem2}, $\phi^{-1}(I) \phi^{-1}(J) \neq \{0\}$. Hence $\{0\} \neq \phi^{-1}(I) \phi^{-1}(J) \subseteq P + Ker(\phi) \subseteq P $. However, $\phi^{-1}(I) \phi^{-1}(J) \subseteq N$ then $\phi^{-1}(I) \phi^{-1}(J) \subseteq P$ since P is a graded weakly prime ideal of N, so $\phi^{-1}(I) \subseteq P$ $\phi^{-1}(J) \subseteq P$. Therefore, $I \subseteq \phi(P)$ or $J \subseteq \phi(P)$. Hence $\phi(P)$ is a graded weakly prime ideal of M. 
\end{theorem}

Next Example and Theorem \ref{thm5}, if $I \subseteq P$ and $\pi : N \rightarrow \Bar{N} := N/I$ is the canonical epimorphism, then $\pi(P)$ is graded weakly prime ideal of $\Bar{N}$ if P is graded weakly prime ideal of N while it is not necessary that P is graded weakly prime ideal in N if $\pi(P)$ is graded weakly prime ideal of $\Bar{N}$.

\begin{example}\label{ex3}
Let $N = \mathbf{Z}_{18}$ be a G-graded near-ring where $N_0 = \mathbf{Z}_{18}$ and $N_1 = \{0\}$. Consider $\pi : N \rightarrow \Bar{N} := N/I$, where $\pi(x)=x$ and $I = \{0, 9\}$. It is easily to check hat $\pi(\{0, 9\}) = \Bar{0}$ is a graded weakly prime ideal of $\Bar{N}$. However, $I \subseteq \{0, 9\}$ is not graded weakly prime ideal of N.
\end{example}

\begin{theorem}\label{thm5}
Let N be a G-graded near-ring and P, I be a graded ideals of N with $I \subseteq P$. Consider $\pi : N \rightarrow \Bar{N} := N/I$ is the canonical epimorphism. If P is a graded weakly prime ideal of N, then $\pi(P)$ is a graded weakly prime ideal of $\Bar{N}$.
\\
\\
\textbf{Proof.} Let J and K be a graded ideals of N with $ \{0\}\neq KJ \subseteq P$, so $\pi(J)$ and $\pi(K)$ are graded ideals of $\Bar{N}$ with $\{0\} \neq \pi(J)\pi(K)=\pi(JK) \subseteq \pi(P)$. Since $\{0\} \neq \pi(J)\pi(K) $ then by Lemma \ref{lem2} $\{0\} \neq \pi^{-1}(\pi(J)) \pi^{-1}(\pi(K))$ and then $\pi^{-1}(\pi(J)) \pi^{-1}(\pi(K)) = JK \subseteq \pi^{-1}(\pi(P))= P+I = P$. Thus $J \subseteq P$ or $K \subseteq P$. Therefore. $J = \pi^{-1}(\pi(J)) \subseteq P \subseteq \pi^{-1}(\pi(P))$ so $\pi(J) \subseteq \pi(P)$ or $\pi(K) \subseteq \pi(P)$. Thus $\pi(P)$ is a grade weakly prime ideal of $\Bar{N}$.
\end{theorem}

Note that, from the definition of graded weakly prime ideals, for any graded ideal of N with $I^2 \subseteq P$ where P is a graded weakly prime ideal of N, if $I \not \subseteq P$ then $I^2 = \{0\}$. If there is another special cases that guarantee $I^2 =\{0\}$. The Theorem \ref{thm6} gives one case of them but before state it the following lemma is presented.

\begin{lemma}\label{lem3}
Let P be a graded weakly prime ideal of N. If $\Bar{I}$ is a graded ideal of N/P with $\Bar{I}\Bar{J} =\{0\}$ for some non zero graded ideal $\Bar{J}$ of N/P. Then there is either $I \subseteq P$ or $PJ = \{0\}$.
\\
\\
\textbf{Proof.} Suppose that $I \not\subseteq P$ and let $p \in P$. Then $(<p>+I)\not\subseteq P$ and $(<p>+I)J \subseteq P$ which implies $(<p>+I)J = \{0\}$ but P is a graded weakly prime ideal of N. Thus $<p>J =\{0\}$ and hence $PJ = \{0\}$.
\end{lemma}

\begin{theorem}\label{thm6}
Let N be a G-graded near-ring and P be a graded weakly prime ideal of N with $P^2 = \{0\}$. If I is a graded ideal of N and $I^2 \subseteq P$, then $I^2= \{0\}$.
\\
\\
\textbf{Proof.} Let P is a graded weakly prime ideal of N and for any $x, y \in I$, there is $<x><y> \subseteq I^2 \subseteq P$. Now, the claim is that $<x><y> = \{0\}$. Suppose not, then since P is graded weakly prime ideal of N, there is $x\in P$ or $y \in P$. If both $x, y \in P$, then $<x><y> \subseteq P^2 =\{0\}$. So, it is assumed that only one of them x or y belongs to P. Take $x \in P$ since $<y><y> \subseteq I^2 \subseteq P$ and by Lemma \ref{lem3} we have $<x><y> \subseteq P<y> = \{0\}$ which implies $I^2 = \{0\}$.
\end{theorem}

Recall that, if $N$ and $M$ is a G-graded near-rings, then $N \times M$ is a G-graded near-ring.  

\begin{theorem}\label{thm7}
Let N and M be a G-graded near-rings and P be a graded ideal of N. Then P is a graded weakly prime ideal of N if and only if $P \times M$ is a graded weakly prime ideal of $N\times M$.
\\
\\
\textbf{Proof.} ($\Rightarrow$) Let P be a graded weakly prime ideal of N and $I \times M$, $J\times M$ be a graded ideal of $N\times M$ such that $\{0\} \neq (I\times M)(J \times M) \subseteq P \times M$. Then $\{0\} \neq (I\times M)(J \times M) = (IJ \times MM) \subseteq P \times M$. So, $\{0\} \neq IJ \subseteq P$ but P is a graded weakly prime ideal of N then $I \subseteq P$ or $J \subseteq P$. Therefore, $I \times M \subseteq P \times M$ or $J \times M \subseteq P \times M$. Thus $P \times M$ is a graded weakly prime ideal of $N \times M$.
\\
($\Leftarrow$) Suppose that $P\times M$ is a graded weakly prime ideal of $N\times M$ and Let I, J be a graded ideals of N such that $\{0\} \neq IJ \subseteq P$. Then $\{0\} \neq (I\times M)(J \times M) \subseteq P \times M$. By assumption we have $I \times M \subseteq P \times M$ or $J \times M \subseteq P \times M$. So $I \subseteq P$ or $J \subseteq P $. Thus P is a graded weakly prime ideal of N.
\end{theorem}

\begin{corollary}\label{coro2}
Let N and M be two G-graded near-rings. If every graded ideal of N and M is a product of graded weakly prime ideals, then every graded ideal of $N \times M$ is a product of graded weakly prime ideals.
\\
\\
\textbf{Proof.} Let I be a graded ideal of N and J be a graded ideal of M such that $I = I_1 ... I_n$ and $J = J_1 ... J_m$ where $I_i$ and $J_j$ is a graded weakly prime ideal of N and M respectively. If the graded ideal is of the form $I \times M$ then $I \times M = (I_1 ... I_n) \times M $ can be written as $(I_1 \times M)...(I_n \times M)$ which is by Theorem \ref{thm7} a product of graded weakly prime ideals. Similarly, if the graded ideal is of the form $N \times J$, then it is a product of graded weakly prime ideals. If the graded ideal is of the form $I \times J$ then it can be written as $(I_1 ... I_n) \times (J_1 ... J_m) = ((I_1 ... I_2) \times M)(N \times (J_1 ... J_m)) = (I_1 \times M)...(I_n \times M)(N \times J_1)...(N \times J_m)$ which is a product of graded weakly prime ideals. 
\end{corollary}

\begin{theorem}\label{thm8}
Let $N$ and $M$ be  two $G$-graded near-rings. Then a graded ideal $P$ of $N \times M$ is graded weakly prime if and only if it has one of the following two forms:
\\
(i) $I \times M$, where $I$ is a graded weakly prime ideal of $N$.
\\
(ii) $N \times J$, where $J$ is a graded  prime ideal of $M$.
\\
\\
\textbf{Proof.} Let $P$ be a graded ideal of $N \times M.$ Then $P$ has one of the following three forms (i) $I \times M$ where $I$ is a graded ideal of $N$ or (ii)$N \times J$, where   $J$ is proper ideal of $M$ or $I \times J$, where  $I\neq N $  and $ J \neq M$. If $P$ is of the form $I \times M$ or of the form $N \times J$
then by Theorem \ref{thm7}, $P$ is graded weakly prime ideal of $N \times N$ if and only if both $I$ and $J$ are graded weakly prime ideals of N and M respectively. Let $P = I \times J$ be a
graded weakly prime ideal of $N \times M$  with  $I \neq N $ and $ J \neq M $. Suppose $ x  \in I $. Then $< x > \times \{0\} \subseteq P$ This implies that either $< x > \times M \subseteq P $ or $(N \times \{0\} \subseteq P$. If $< x > \times M \subseteq P $, then
 $M = J$ and if $(N \times \{0\} \subseteq P$, then $N = I$ this is a contradiction. Hence $I\times J$ can not be graded weakly prime ideal of $N \times M$ if both $I$ and $J$ are graded ideals.
\end{theorem}

\begin{theorem}\label{thm9}
Let $N$ and $M$ be two G-graded near-rings. Then $P= \{0\} \times \{0\}$ is a graded weakly prime ideal of $N \times M$.
\\
\\
\textbf{Proof.} Let $I = \{0\}$ be a graded ideal of N. Suppose that $x \in I-\{0\}$ then $<x> \times \{0\} \subseteq P$ and $<x> \times \{0\} \neq \{0\}$. This implies that either $<x> \times M \subseteq P$ or $N \times \{0\} \subseteq P$ then $N = I$ this is a contradiction. So $I -\{0\}$ is empty. Similarly if $J = \{0\}$ where J is a graded ideal of M. Therefore, P is a graded weakly prime ideal of $N \times M$.
\end{theorem}

\begin{proposition}\label{prop5}
Let N be a G-graded near-ring and P, I be two graded ideals of N. If P and I are graded weakly prime ideal of N, then $P \cup I$ is a graded weakly prime ideal of N.
\\
\\
\textbf{Proof.} Let J and K be a two graded ideals of N such that $\{0\}\ \neq JK \subseteq P \cup I$. Then $JK \subseteq P$ or $JK \subseteq I$. Since $JK \neq \{0\}$ and if $JK \subseteq P$, then $J \subseteq P$ or $K \subseteq P$ since P is a graded weakly prime ideal of N. Hence $J \subseteq P \cup I$ or $K \subseteq P \cup I$. If $JK \subseteq I$, then $J \subseteq I$ or $K \subseteq I$ since I is a graded weakly prime ideal of N. thus $J \subseteq P \cup I$ or $K \subseteq P\cup I$. Therefore, $P\cup I$ is a graded weakly prime ideal of N.
\end{proposition}

\section{GRADED ALMOST PRIME IDEALS OVER GRADED NEAR RINGS}

In this section, we introduce graded almost prime ideals over near-rings concept and study their basic properties.

\begin{definition}\label{def2}
Let G be a multiplucative monoid group with identity element and N be a G-graded near-ring. A graded ideal P of N is called graded almost prime ideal of N if whenever $IJ\subseteq P$ and $IJ \not\subseteq (P^2 \cap N)$, then either $I \subseteq P$ or $J \subseteq P$, for any graded ideals I and J in N.
\end{definition}

\begin{example}\label{ex4}
Consider a G-graded near-ring which is defined in Example \ref{ex1}. Not that $P_4= \{0,4,8\}$ is a graded almost prime ideal of N but not graded weakly prime ideal of N. However, $P_5 = \{0,6\}$ is neither graded weakly prime ideal of N nor graded almost prime ideal of N. 
\end{example}

In the previous section, it was observed that if P is graded prime ideal then it is graded weakly prime ideal but the converse is not true for example P1 in Example \ref{ex1} is not graded prime ideal of N but it is graded weakly prime ideal of N. Also, by Example \ref{ex4}, a graded almost prime ideal of N may not implies a graded weakly prime ideal of N . Now, the question is: Does a graded weakly prime ideal give a graded almost prime ideal? The next theorem answers this question.

\begin{theorem}\label{thm10}
Let N be a G-graded near-ring and P be a graded ideal of N. If P is a graded weakly prime ideal of N, then P is a graded almost prime ideal of N.
\\
\\
\textbf{Proof.} Let P be a graded weakly prime ideal of N and I, J be two graded ideals of N such that $IJ\subseteq P$ and $IJ \not\subseteq (P^2 \cap N)$. If P is a graded prime ideal of N then P is a graded almost prime ideal of N. Otherwise, $(P^2 \cap N) = \{0\}$ by Theorem \ref{thm1} $IJ \neq \{0\}$ since $IJ \not\subseteq (P^2 \cap N) = \{0\}$. But P is a graded weakly prime ideal of N. Therefore, either $I\subseteq P$ or $J \subseteq P$ which means that P is a graded almost prime ideal of N.  
\end{theorem}

\begin{proposition}\label{prop6}
Let N be a G-graded near-ring and P be a graded prime ideal of N. If P is a graded almost prime ideal of N and $((P^2 \cap N):P) \subseteq P$, then P is a graded prime ideal of N.
\\
\\
\textbf{Proof.} Suppose that P is not graded prime ideal of N. Then there exist $I \not\subseteq P$ and $J \not\subseteq P$ satisfying that $IJ \subseteq P$ where I and J are two graded ideals of N. If $IJ \not\subseteq (P^2 \cap N)$ we are done. So, it is assumed $IJ \subseteq (P^2 \cap N)$. Consider $I(J+P) \subseteq P$ if $I(J+P) \not\subseteq (P^2 \cap N)$ then there is  $I \subseteq P$ or $J \subseteq P$ this is a contradiction. Otherwise, $I(J+P) \subseteq (P^2 \cap N)$ then $IP \subseteq (P^2 \cap N)$ which implies $I \subseteq ((P^2 \cap N):P) \subseteq P$ which is a contradiction. Thus P is a graded prime ideal of N.
\end{proposition}

Next, some equivalent conditions are given for a graded ideal to be graded almost prime ideal in the G-graded near-ring.

\begin{theorem}\label{thm11}
Let N be a G-graded near-ring and P be a graded ideal of N. Then the following are equivalent:
\\
(1) For $x,y \ and\ z \in N$ with $ x(<y>+<z>) \subseteq P$ and $x(<y>+<z>) \not \subseteq (P^2\cap N)$ there is $x\in P$ or $y \ and\ z \in P$.
\\
(2) For $x\in N$ but $x \not \in P$, $(P:<x>+<y>) = P \cup ( (P^2 \cap N):<x>+<y>)$ for any $y \in N$.
\\
(3) For $x\in N$ but $x \not\in P$ we have $(P:<x>+<y>) = P$ or $(P:<x>+<y>) = ((P^2\cap N):<x>+<y>)$ for any $y\in N$.
\\
(4) P is a graded almost prime ideal of N.
\\
\\
\textbf{Proof.} $(1) \Rightarrow (2)$: Let $t \in N$ and $t \in (P :< x > + < y >)$ for any y and x belongs to N but $x \not\in P$. Then $t(< x > + < y >) \subseteq P$. If $t(< x > + < y >) \subseteq (P^2\cap N)$. Then $t \in ((P^2 \cap N) :< x > + < y >)$. Otherwise, we get $ t(< x > + < y >) \not\subseteq (P^2\cap N)$. Thus, $t \in P$ by hypothesis. 
\\
$(2) \Rightarrow (3)$: It is following directly from Lemma \ref{lem1}.
\\
$(3) \Rightarrow (4)$: Let I and J be a graded ideal of N such that $IJ \subseteq P$ and $IJ \not\subseteq (P^2\cap N)$. Suppose that $I \not\subseteq P$ and $J \not\subseteq P$. Then there exist $j \in J$ with $j \not\in P$. Now, it is claimed that $IJ \subseteq (P^2 \cap N)$. Let $j_1 \in J$, then $I(<j>+<j_1> \subseteq P$, which implies $I \subseteq (P:<j>+<j_1>)$. Then by assumption, $I(<j>+<j_1>) \subseteq (P^2\cap N) $ which gives $Ij_1 \subseteq (P^2 \cap N) $. Thus $IJ \subseteq (P^2\cap N)$ and hence P is a graded almost prime ideal of N.
\\
$(4) \Rightarrow (1)$: If $ x(<y>+<z>) \subseteq P$ and $ x(<y>+<z>) \not\subseteq (P^2\cap N)$, then $\{0\} \neq <x>(<y>+<z>) \subseteq P$. Since P is a graded almost prime ideal of N, there is $<x> \subseteq P$ or $<y>+<z> \subseteq P$. By assumption, $x, y \ and\ z \in N$. Hence $x\in P$ or $y \ and\ z \in P$. 
\end{theorem}

\begin{theorem}\label{thm12}
Let N be a G-graded near-ring and P be a graded ideal of N. Then the following are equivalent:
\\
(1) P is a graded almost prime ideal of N.
\\
(2) For any ideals I and J in N with $P \subset I$ and $P \subset J$, then there is either $IJ \subseteq (P^2\cap N)$ or $IJ \not\subseteq P$.
\\
(3) For any ideals I and J in N with $I \not\subseteq P$ and $J \not\subseteq P$, then there is either $IJ \subseteq (P^2\cap N)$ or $IJ \not\subseteq P$.
\\
\\
\textbf{Proof.}$(1) \Rightarrow (2)$: Let I an J be two graded ideals of N with $P \subset I$, $P \subset J$ and $IJ\not\subseteq (P^2\cap N)$. Take $i \in I$ and $j \in J$ with $i \not \in P$ and $j \not \in P$, which implies that $<i><j> \not\subseteq P$ and hence $IJ \not \subseteq (P^2\cap N)$ and $IJ \not\subseteq P$. 
\\
$(2) \Rightarrow (3)$: Let I and J be a graded ideals of N with $I \not\subseteq P$ and $J \not\subseteq P$. Then there exists $i_1 \in I$ and $j_1 \in J$ such that $i_1 \not \in P$ and $j_1 \not\in P$. Suppose that $<i><j> \not (P^2\cap N)$ for some $i \in I$ and $j \in J$. Then $(P+<i>+<i_1>)(P+<j>+<j_1>) \not \subset (P^2 \cap N)$ and $P \subset (P+<i>+<i_1>) $, $ P \subset (P+<j>+<j_1>)$. By hypothesis, $(P+<i>+<i_1>)(P+<j>+<j_1>) \not\subseteq P$. So, $<i>(P+<j>+<j_1>) +<i_1>(P+<j>+<j_1>) \not\subseteq P$. Hence there exists $i' \in <i>$, $i_1'\in <i_1>$, $j', j'' \in <j>$, $j_1', j_1'' \in <j_1>$ and $p_1, p_2 \in P$ such that $i'(p_1 + j' + j_1') + i_1'(p_2 + j'' + j_1'') \not\in P$. Thus  $i'(p_1 + j' + j_1') - i'(j' + j_1') + i'(j' + j_1') + i_1'(p_2 + j'' + j_1'') - i_1'(j'' + j_1'') + i_1'(j'' + j_1'') \not\in P$. But $i'(p_1 + j' + j_1') - i'(j' + j_1') \in P$ and $i_1'(p_2 + j'' + j_1'') - i_1'(j'' + j_1'') \in P$. Which implies neither $i'(j' + j_1')$ nor $i_1'(j'' + j_1'')$ belongs to P. Therefore, $IJ \not\subseteq P$.
\\
$(3) \Rightarrow (1) $: Follows directly from the definition of graded almost prime ideals of N.
\end{theorem}

\begin{proposition}\label{prop7}
Let N be a G-graded near-ring, A be a totally ordered set and $(P_a)_{a\in A}$ be a family of graded almost prime ideals of N with $P_a \subseteq P_b$ for any $a,b \in A$ with $a\leq b$. Then $P = \bigcap\limits_{a\in A} P_a$ is a graded almost prime ideal of N.
\\
\\
\textbf{Proof.} Let I and J be two graded ideals of N with $ IJ \subseteq P$ but $IJ \not\subseteq (P^2\cap N)$, which implies for all $a \in A$ there is $IJ \subseteq P_a$. If there exists $a \in A$ such that $I \not\subseteq P_a$, then $J \subseteq P_a$. Hence for all $a \leq b$ there is $J \subseteq P_b$. If there exists $c < a$ such that $j \not \subseteq P_c$, then $I \subseteq P_c $ and then $I \subseteq P_a$, this is a contradiction. Hence for any $a \in A$, there is $J \subseteq P_a$. Therefore, $J \subseteq P$. 
\end{proposition}

\begin{proposition}\label{prop8}
Let N be a G-graded near-ring and P be an intersection of some graded almost prime ideals of N. Then for any graded ideal I of N satisfying that $ I^2 \subseteq P$  but $I^2 \not\subseteq (P^2\cap N)$ we have $I \subseteq P$.
\\
\\
\textbf{Proof.} Let $P_a$ be a set of graded almost prime ideals of N, P be the intersection of $P_a$ and I be a graded ideal of N such that $ I^2 \subseteq P$ but $I^2 \not\subseteq (P^2\cap N)$. Then $I^2$ is subset of each $P_a$ since $P_a$ is graded almost prime ideal of N there is $I \subseteq P_a$. Therefore, $I \subseteq P$.
\end{proposition}

\begin{lemma}\label{lem4}
Let N and M be two G-graded near-rings and $\phi$ be a surjective homomorphism from N into M. For any two graded ideals P, I and J of N if $IJ \not\subseteq P$, then $\phi^{-1}(I) \phi^{-1}(J) \not\subseteq \phi^{-1}(P)$.
\\
\\
\textbf{Proof.} If $I \not\subseteq P$, then $\phi^{-1}(I) \not\subseteq \phi^{-1}(P)$ since $\phi(I) \subseteq \phi(P)$, then $I = \phi(\phi^{-1}(I)) \subseteq \phi(\phi^{-1}(P)) \subseteq P$. Hence if $IJ \not\subseteq P$, then $\phi^{-1}(IJ)=\phi^{-1}(I) \phi^{-1}(J) \not \subseteq \phi^{-1}(P)$. 
\end{lemma}

\begin{theorem}\label{thm13}
Let N and M be two G-graded near-rings and $\phi$ be a surjective homomorphism from N into M. Then the image of graded almost prime ideal of N which contains the kernal of $\phi$ is a graded almost prime ideal of M.
\\
\\
\textbf{Proof.} Suppose that $IJ \subseteq \phi(P)$ and $IJ \not \subseteq ((\phi(P))^2 \cap N)$ where I and J be two graded ideals of N and P is a graded almost prime ideal of N. By Lemma \ref{lem4} $\phi^{-1}(I) \phi^{-1}(J) \not \subseteq (P^2\cap N)$. Hence $\phi^{-1}(I) \phi^{-1}(J)  \subseteq P+Ker(\phi) = P$. Since P is a graded almost prime ideal of N then $\phi^{-1}(I) \subseteq P$ or $ \phi^{-1}(J) \subseteq P$. Therefore, $I \subseteq \phi(P)$ or $J \subseteq \phi(P)$. Hence $\phi(P)$ is a graded almost prime ideal of M.  
\end{theorem}

\begin{theorem}\label{thm14}
Let N be a G-graded near-ring and P, I be a graded ideals of N with $I \subseteq P$. Consider $\pi : N \rightarrow \Bar{N} := N/I$ is the canonical epimorphism. If P is a graded almost prime ideal of N, then $\pi(P)$ is a graded almost prime ideal of $\Bar{N}$.
\\
\\
\textbf{Proof.} Let P be a graded almost prime ideal of N and $\pi(J), \pi(K)$ be a graded ideal of $\Bar{N}$ with $\pi(J)\pi(K)= \pi(JK) \subseteq \pi(P)$ and $\pi(J)\pi(K)= \pi(JK) \not\subseteq ((\pi(P))^2 \cap N)$. Since $\pi(J)\pi(K)= \pi(JK) \not\subseteq ((\pi(P))^2 \cap N)$ then by Lemma \ref{lem4} $JK \not \subseteq P^2$. Therefore, $\pi^{-1}(\pi(J)) \pi^{-1}(\pi(K)) = JK \subseteq \pi^{-1}(\pi(P))= P+I = P$. Hence $J \subseteq P$ or $K \subseteq P$ and hence $\pi(J) \subseteq of \pi(P)$ or $\pi(K) \subseteq of \pi(P)$. Thus $\pi(P)$ is a graded almost prime ideal of $\Bar{N}$.
\end{theorem}

\begin{lemma}\label{lem5}
Let P be a graded almost prime ideal of N. If $\Bar{I}$ is a graded ideal of N/P with $\Bar{I}\Bar{J} =\{0\}$ for some non zero graded ideal $\Bar{J}$ of N/P. Then there is either $I \subseteq P$ or $PJ \subseteq (P^2 \cap N)$.
\\
\\
\textbf{Proof.} Suppose that $I \not\subseteq P$ and let $p \in P$. Then $(<p>+I)\not\subseteq P$ and $(<p>+I)J \subseteq P$ which implies $(<p>+I)J \subseteq (P^2\cap N)$ but P is a graded almost prime ideal of N. Thus $<p>J \subseteq (P^2 \cap N)$ and hence $PJ \subseteq (P^2 \cap N)$.
\end{lemma}

\begin{theorem}\label{thm15}
Let N be a G-graded near-ring and P be a graded almost prime ideal of N. If I is a graded ideal of N and $I^2 \subseteq P$, then $I^2 \subseteq (P^2\cap N)$.
\\
\\
\textbf{Proof.} Let P is a graded almost prime ideal of N and for any $x, y \in I$ we have $<x><y> \subseteq I^2 \subseteq P$. Now, the claim is that $<x><y> \subseteq (P^2\cap N)$. Suppose not, then since P is graded almost prime ideal of N, there is $x\in P$ or $y \in P$. So, it is assumed that only one of them x or y belongs to P. Take $x \in P$ since $<y><y> \subseteq I^2 \subseteq P$ and by Lemma \ref{lem5} we have $<x><y> \subseteq P<y> \subseteq (P^2\cap N)$ which implies $I^2 \subseteq (P^2 \cap N)$.
\end{theorem}

Previous theorem is important in some G-graded near-ring with unique maximal in $N$ (as $N$ is a G-graded near-ring). If this maximal ideal is graded ideal and satisfying that $M M = M^2 \cap N$ like the G-graded near-ring $N= (\mathbf{Z}_{16}, + , \times)$ with $G$ defined as Example \ref{ex1} where $N_0 = \mathbf{Z}_{16}$ and $N_1 = \{0\}$. Note that $N$ as a G-graded ring has unique Maximal ideal $M=\{0, 2, 4, 6, 8, 10, 12, 14\}$ and $M$ satisfies the property $M M = (M^2 \cap N)$. The importance of such G-graded near-rings is explained in the following theorem.

\begin{theorem}\label{thm16}
Let N be a G-graded near-ring. If M is the unique maximal ideal of N with $M M = (M^2 \cap N)$, then for any graded ideal P of N with $M^2\cap N \subseteq P$. We have P is a graded almost prime ideal of N if and only if $(M^2 \cap N) = (P^2 \cap N)$.
\\
\\
\textbf{Proof.} ($\Rightarrow$) Let P be a graded almost prime ideal of N. Then there is $M^2 \cap N = M M \subseteq (P^2\cap N)$ by Theorem \ref{thm15}. Also, $(P^2 \cap N) = (M^2 \cap N)$ since M is the unique Maximal ideal of N.
\\
($\Leftarrow$) Let $(M^2 \cap N) = (P^2 \cap N)$ we claim that P is a graded almost prime ideal of N. Let I and J be two graded ideals of N with $IJ \subseteq P$. Since M is the unique Maximal ideal of N so $I \subseteq M$ and $J \subseteq M$. Therefore, $IJ \subseteq M M = (M^2 \cap N) = (P^2 \cap N)$. Hence P is a graded almost prime ideal of N.
\end{theorem}

\begin{theorem}\label{thm17}
Let N and M be a G-graded near-rings and P be a graded ideal of N. Then P is a graded almost prime ideal of N if and only if $P \times M$ is a graded almost prime ideal of $N\times M$.
\\
\\
\textbf{Proof.} ($\Rightarrow$) Let P be a graded almost prime ideal of N and $I \times M$, $J\times M$ be a graded ideal of $N\times M$ such that $ (I\times M)(J \times M) \subseteq P \times M$ and $ (I\times M)(J \times M) \not\subseteq ((P \times M)^2 \cap (N\times M)$ . Then $(I\times M)(J \times M) = (IJ \times MM) \subseteq P \times M$ and $(I\times M)(J \times M) = (IJ \times MM)\not\subseteq ((P^2  \cap N)\times M)$. So, $ IJ \subseteq P$ and $IJ \not \subseteq (P^2\cap N)$  but P is a graded almost prime ideal of N then $I \subseteq P$ or $J \subseteq P$. Therefore, $I \times M \subseteq P \times M$ or $J \times M \subseteq P \times M$. Thus $P \times M$ is a graded almost prime ideal of $N \times M$.
\\
($\Leftarrow$) Suppose that $P\times M$ is a graded almost prime ideal of $N\times M$ and Let I, J be a graded ideals of N such that $ IJ \subseteq P$ and $IJ \not\subseteq (P^2\cap N)$. Then $ (I\times M)(J \times M) \subseteq P \times M$ and $ (I\times M)(J \times M) \not\subseteq ((P \times M)^2 \cap (N\times M)$. By assumption we have $I \times M \subseteq P \times M$ or $J \times M \subseteq P \times M$. So $I \subseteq P$ or $J \subseteq P $. Thus P is a graded almost prime ideal of N.
\end{theorem}

\begin{corollary}\label{coro3}
Let N and M be two G-graded near-rings. If every graded ideal of N and M is a product of graded almost prime ideals, then every graded ideal of $N \times M$ is a product of graded almost prime ideals.
\\
\\
\textbf{Proof.} Let I be a graded ideal of N and J be a graded ideal of M such that $I = I_1 ... I_n$ and $J = J_1 ... J_m$ where $I_i$ and $J_j$ is a graded almost prime ideal of N and M respectively. If the graded ideal is of the form $I \times M$ then $I \times M = (I_1 ... I_n) \times M $ can be written as $(I_1 \times M)...(I_n \times M)$ which is by Theorem \ref{thm17} a product of graded almost prime ideals. Similarly, if the graded ideal is of the form $N \times J$, then it is a product of graded almost prime ideals. If the graded ideal is of the form $I \times J$ then it can be written as $(I_1 ... I_n) \times (J_1 ... J_m) = ((I_1 ... I_2) \times M)(N \times (J_1 ... J_m)) = (I_1 \times M)...(I_n \times M)(N \times J_1)...(N \times J_m)$ which is a product of graded almost prime ideals. 
\end{corollary}

\bibliographystyle{amsplain}

\end{document}